\documentclass[11pt,twoside,center]{article}

\usepackage{amsmath}%
\usepackage{amsfonts}%
\usepackage{amssymb}%
\usepackage{graphicx}%
\usepackage{listings}
\usepackage{anysize}
\usepackage{float}
\restylefloat{figure}
\usepackage{graphicx}

\usepackage[left=5.0cm, right=4.0cm, top=4.0cm, bottom=5.0cm]{geometry}

\newtheorem{definition}{Definition}
\newtheorem{theorem}{Theorem}
\newtheorem{example}{Example}
\newtheorem{lemma}{Lemma}

\newtheorem{remark}{Remark}
\newtheorem{proposition}{Proposition}
\newtheorem{algorithm}{Algorithm}
\newenvironment{proof}[1][Proof]{\textbf{#1.} }

\title{\bf A Hamiltonian approach to implicit systems, generalized solutions and applications}

\author{Dan TIBA\footnote{
{\tt dan.tiba@imar.ro} Academy of Romanian Scientists and Institute of Mathematics, Romanian Academy, P.O. BOX 1-764, 014700 Bucharest, Romania}}
\date{~}
\begin{document}

\maketitle

\pagestyle{myheadings}
\markboth{Dan Tiba
}{\textbf{A Hamiltonian approach}}

\bigskip

\begin{abstract}
We introduce a constructive method that provides the local solution of general implicit systems in arbitrary dimension via Hamiltonian type equations. A variant of this approach constructs  parametrizations of the manifold, extending the usual implicit functions solution.  We also investigate the unsolved critical case of the implicit functions theorem, define the notion of \emph{generalized solution} and prove existence and basic properties. Relevant examples and counterexamples are also indicated. The applications concern new necessary conditions (with less Lagrange multipliers), perturbations and algorithms in non convex optimization problems.
\newline
\end{abstract}

{\bf keywords:}  local parametrizations, uniqueness, critical case, generalized solutions, nonlinear programming.

\section{Introduction}
\label{intro}
In the Euclidean space $R^d$, $d \in N$, we consider the general implicit functions system:

\begin{center}
\begin{eqnarray}\label{eq11}
F_1(x_1, \ldots, x_d) = 0,\nonumber\\
F_2(x_1, \ldots, x_d) = 0,\\
\ldots\ldots\ldots \ldots \ldots \ldots\nonumber\\
F_l(x_1, \ldots, x_d) = 0,\nonumber
\end{eqnarray}
\end{center}

\noindent where $l \in N$, $1 \leq l \leq d-1$ and $F_j \in C^1(\Omega)$, $F_j(x^0) = 0, \; j = \overline{1, l}$, $x^0 \in \Omega \subset R^d$ bounded domain, given.

The problem (\ref{eq11}) has a long and well known history and we quote the monographs of Krantz and Parks \cite{c6}, Dontchev and Rockafellar \cite{c3} for a comprehensive presentation, including important applications and recent research developments. We also mention the book by Thorpe \cite{c16}, where related ideas are discussed from the point of view of differential geometry. In particular, it is known that one can associate to (\ref{eq11}) a system of nonlinear (partial) differential equations (basically derived from the differentiation formula, under usual assumptions), see \cite{c6}, Ch. 4.1.

In the recent paper \cite{c17}, in dimension two and three, it was shown that one can associate to (\ref{eq11}) other (essentially simpler) systems of ordinary differential equations, under the mere assumption that $F_j \in C^1(\Omega)$, $j = \overline{1, l}$ and in the absence of any independence-type condition. These new systems provide a  constructive (local) parametrization of the solution of (\ref{eq11}) around $x^0$ under certain conditions.

Moreover, it is possible to define a local generalized solution of (\ref{eq11}) even in the critical case, in arbitrary dimension. Our approach to this old question is novel. A variant  (in dimension three) was discussed in \cite{c11} as well, where it was proved that it is enough to use ordinary differential systems (especially of Hamiltonian type) in order to solve locally (\ref{eq11}) via appropriate parametrizations of the unknowns. Several relevant numerical examples are also indicated.

In this paper, we discuss the solution of the general implicit system (\ref{eq11}) in arbitrary dimension, by using a new iterated system of ordinary differential equations. The approach has a constructive character and we indicate two variants that give a parametrization of the unknowns in (\ref{eq11}) or construct exactly the classical implicit solution (in function form, see Theorem \ref{th27}). Obtaining parametrizations is advantageous since they may provide a better description of the manifold. This is done in Section 2, under the usual nondegeneracy condition from the implicit function theorem. We underline that the existence question is well known (via the classical implicit functions theorem or the inverse function theorem, etc.), but a general and effective construction seems not to be available, to the best of our knowledge. The systems of ordinary differential equations that we use here  are derived from a first order partial differential system of equations. An interesting fact is that, although just continuity is valid for the right-hand side, we also prove the uniqueness of the solution. This may be compared with the results from \cite{bd}, \cite{dpl} (see Remarks \ref{r1}, \ref{r2}).  As a first application, in the final part of Section 2, perturbations of the system
(\ref{eq11}) are investigated, both for implicit functions and implicit parametrizations.

\noindent
We also recall that, in algebraic geometry, implicitization and parametrization (via rational functions) are important subjects, Gao \cite{c4}, Wang \cite{c18}, Schicho \cite{c15}. General parametrization methods are not known, Gao \cite{c4} and recent papers study approximate parametrization approaches, Dobiasova \cite{c2}, Yang, J\" uttler and Gonzales-Vega \cite{c19}.

In Section 3, we show how to solve the critical case as well, under the assumption $F_j \in C^1(\Omega)$, $j = \overline{1, l}$, in the absence of (\ref{eq21}). We introduce the notion of (local) \emph{generalized solution}, prove its existence and basic properties. We also indicate some relevant examples. The generalized solution obtained by our method covers all possible cases and is an extension of the notion of  local solution from the implicit function theorem, in the classical case. Singular situations in the implicit functions theorem were discussed by different methods in \cite{a}, \cite{b}, \cite{c} and a comprehensive account can be found in \cite{c6}, Ch.5.4, where it is specified that a complete solution of the critical case is not known.

The implicit parametrizations constructed in this work are also useful in computations of integrals on implicitly defined  manifolds, which seems to be yet unsolved. There is a recent interest in this direction due to important questions in the well known level set method for evolving surfaces, see \cite{kt} and its references, or in shape optimization, \cite{ti}.

The last section is devoted to applications in nonlinear programming. We use reduced gradients to obtain optimality conditions in the simpler Fermat form, involving no or fewer multipliers, i.e. the constraints can be eliminated, at least partially.  We also introduce two algorithms, including one that works in the critical case as well. Some numerical examples are also provided together with a comparison with the relaxation approach \cite{ssb} or the fmincon routines in MatLab.

\label{sec:1}

\section{Implicit  parametrizations}

In this section, we discuss the system (\ref{eq11}) under the classical independence assumption. To fix ideas, we assume
\begin{equation}\label{eq21}
\displaystyle\frac{D(F_1, F_2, \ldots, F_l)}{D(x_1, x_2, \ldots, x_l)} \neq 0 \quad {\rm in} \; x^0 = (x_1^0, x_2^0, \ldots, x_d^0). 
\end{equation}

The hypothesis (\ref{eq21}) will be dropped in the next section.

Clearly, condition (\ref{eq21}) remains valid on a neighbourhood $V \in \mathcal{V}(x^0)$, $V \subset \Omega$, under the $C^1(\Omega)$ assumption on $F_j(\cdot)$, $j = \overline{1, l}$ and we denote by $A(x), x  \in V$, the corresponding nonsingular $l \times l$ matrix from (\ref{eq21}).

We introduce on $V$ the undetermined linear systems of equations with unknowns $v(x) \in R^d, \; x \in V$:

\begin{equation}\label{eq22}
v(x) \cdot \nabla F_j(x) = 0, \quad j = \overline{1, l}. 
\end{equation}

We shall use $d-l$ solutions of (\ref{eq22}) obtained by fixing successively the last $d-l$ components of the vector $v(x) \in R^d$ to be the rows of the identity matrix in $R^{d-l}$ multiplied by $\Delta (x) = {\rm det} A(x)$. Then, the first $l$ components are uniquely determined, by inverting $A(x)$, due to (\ref{eq21}).

In this way, the obtained $d-l$ solutions of (\ref{eq22}), denoted by $v_1(x), \ldots$, $v_{d-l}(x) \in R^d$, are linear independent, for any $x \in V$.

Moreover, these vector fields are continuous in $V$ as $\nabla F_j(\cdot)$ are continuous in $V$ and the Cramer's rule ensures the continuity of the solution for linear systems with respect to the coefficients. Other choices of solutions for (\ref{eq22}), useful in this section, are possible  (see Theorem \ref{th27}).

We introduce now $d-l$ nonlinear systems of first order partial differential equations associated to the vector fields $(v_j(x))_{j = \overline{1, d-l}}$, $x \in V \subset \Omega$. Furthermore, we denote the sequence of independent variables by $t_1, t_2, \ldots, t_{d-l}$.

These systems have a nonstandard (iterated) character in the sense that the solution of one of them is used as initial condition in the next one. Consequently, the independent variables in the "previous" systems enter as parameters in the next system via the initial conditions. Due to their simple structure (just one derivative in each equation), we stress that each system (\ref{eq23}), (\ref{eq24}),..., (\ref{eq25}), may be interpreted as an ordinary differential system in $ V \subset R^d $, with parameters, although partial differential notations are used:

\begin{equation}\label{eq23}
\displaystyle\frac{\partial y_1(t_1)}{\partial t_1} = v_1(y_1(t_1)), \quad t_1 \in I_1 \subset R, \end{equation}

$$
y_1(0) = x^0;
$$

\begin{equation}\label{eq24}
\displaystyle\frac{\partial y_2(t_1, t_2)}{\partial t_2} = v_2(y_2(t_1, t_2)), \quad t_2 \in I_2(t_1) \subset R, 
\end{equation}

$$
y_2(t_1, 0) = y_1(t_1);
$$
$$
\ldots \quad \ldots \quad \ldots \quad \ldots \quad \ldots \quad \ldots \quad \ldots \quad \ldots \quad \ldots
$$

\begin{equation}\label{eq25}
\displaystyle\frac{\partial y_{d-l}(t_1, t_2, \ldots, t_{d-l})}{\partial t_{d-l}} = v_{d-l}(y_{d-l}(t_1, t_2, \ldots, t_{d-l})), 
\end{equation}
$$
t_{d-l} \in I_{d-l}(t_1, \ldots, t_{d-l-1}),
$$
$$
y_{d-l}(t_1, \ldots, t_{d-l-1}, 0) = y_{d-l-1}(t_1, t_2, \ldots, t_{d-l-1}).
$$

\vspace{3mm}
Here, the notations $I_1, I_2(t_1),. \ldots, I_{d-l}(t_1, \ldots, t_{d-l-1})$ are $d-l$ real intervals, containing 0 in interior and depending, in principle, on the "previous" parameters. The existence of the solutions $y_1, y_2, \ldots, y_{d-l}$ follows by the Peano theorem due to the continuity of the vector fields $(v_j)_{j = \overline{1, d-l}}$ on $V$. We show in the next two results that all these subsystems have the uniqueness property as well.

\begin{theorem}\label{th1}
Under assumption (\ref{eq21}), if $l = d-1$, then the system (\ref{eq23})-(\ref{eq25}) consists just of one subsystem of dimension $d$ with the uniqueness property.
\end{theorem}

\begin{proof}
By the implicit functions theorem, under hypothesis (\ref{eq21}), around $x^0 = (x_1^0, x_2^0, \ldots, x_d^0)$, there are some open sets (expressed as products of open intervals) such that
$x_1 = \tilde{x}_1(x_d), \dots,  x_{d-1} = \tilde{x}_{d-1}(x_d)$, uniquely defined. Replacing these expressions in the last equation of the differential system, we get an ordinary differential equation in $x_d$ with continuous right-hand side, that has a unique solution (it can be integrated in fact) on some interval. The essential point here is that the right-hand side is non zero around the origin, due to (\ref{eq21}) and one can divide by it, etc. Next, the above equalities give all the components, uniquely determined on some interval around the origin.
\end{proof}

\begin{remark}\label{r1} 
\end{remark}
\vspace{-3mm}
The cases $d=1$ (one ODE), $ d=2, l=1$ (one Hamiltonian system) are discussed in \cite{bd}, using similar ideas.

\begin{theorem}\label{th2}
Under assumption (\ref{eq21}), if $1 \leq l \leq d-2, d \geq 3$, every subsystem in (\ref{eq23})-(\ref{eq25}) has a unique solution.

\end{theorem}

\begin{proof}
We shall proceed by induction on $d$, while $l \leq d-2$ is arbitrarily fixed. In Theorem \ref{th1} and Remark \ref{r1}, we have clarified the small dimension cases, including $d=3, l=2$. Here, we start the induction with $d = 3, l = 1$, which also gives a hint on the general argument.

Denoting by $(x_0, y_0, z_0)$ the initial condition in this case, using (\ref{eq22}) and (\ref{eq23})-(\ref{eq25}) and fixing that $F_x(x_0, y_0, z_0) \neq 0$ in (\ref{eq21}), we obtain two iterated Hamiltonian systems (see \cite{c11}, \cite{c17}). We write just the first one (the second one is similar):

\vspace{-2mm}
\begin{eqnarray}
x^{\prime}&=&-F_y(x,y,z),\hspace{1.95cm}t\,\in I_1,\nonumber\\
y^{\prime}&=&F_x(x,y,z),\hspace{2.2cm}t\,\in I_1,\,\,\,\nonumber\\
z^{\prime}&=&0,\hspace{3.5cm}t\in I_1,\nonumber\\
x(0)&=&x_0,\,\,y(0)\,\,\,\,=\,\,\,\,y_0,\,\,z(0)\,\,\,\,=\,\,\,\,z_0;\label{eq7}
\end{eqnarray}

It is clear that $z(t) = z_0$ is the constant solution and what remains is in fact a simple Hamiltonian system with $d=2, l=1$ that has unique solution by Theorem \ref{th1}.

\noindent
In the general case, we assume that the statement is valid for $d-1$ and any arbitrarily given $l \leq d-3$.  One can reduce any subsystem in (\ref{eq23})-(\ref{eq25}) from dimension $d$ to $d-1$ as above in (\ref{eq7}) since the construction of the right-hand  side $v_j(x)$ involves $d-1-l \geq 1$ null components. If we restrict the implicit system (\ref{eq11}) to $R^{d-1}$ by fixing the existing constant independent variable, then the hypothesis (\ref{eq21}) remains valid for the reduced implicit system and $A(x)$ is again the nonsingular matrix giving the maximal rank. We notice that the modified  $v_j(x)$ satisfies as well the reduced algebraic system (\ref{eq22}) since we have just removed one null component of the original $v_j(x)$. Therefore, the reduced differential subsystem comes as well from a reduced implicit system satisfying (\ref{eq21}).

\noindent
We have the following two variants in dimension $d$. The first one  is $l=d-2$. Then, in dimension $d-1$, the uniqueness of the reduced differential system follows by Theorem \ref{th1}. If $l \leq d-3$, then we can apply the hypothesis of the induction and again we get uniqueness.
This ends the proof.

\end{proof}

\begin{remark}\label{r2} 
\end{remark}
\vspace{-3mm}
Note that for systems associated to divergence free fields, the uniqueness results of \cite{dpl} are valid under certain Sobolev type regularity conditions. However, under our hypotheses, we have just continuity in the right-hand side of the differential system (\ref{eq23})-(\ref{eq25}) and \cite{dpl} cannot be applied.

\vspace{2mm}
\noindent
The next theorems prove more properties and clarify the use of the above setting. 

\begin{theorem}\label{th21}
a) There are closed intervals $I_j \subset R$, $0 \in {\rm int} I_j$, independent of the parameters, such that $I_j \subset I_j(t_1, t_2, \ldots, t_{j-1})$, $ j = \overline{1, d-l}$.

\vspace{1mm}
b) The solutions of the systems (\ref{eq23}) - (\ref{eq25}) are of class $C^1$ in any existence point and we have:
\vspace{2mm}

$\displaystyle\frac{\partial y_{d-l}} {\partial t_k}(t_1, \ldots ,t_{d-l}) =  v_k(y_{d-l}(t_1, \ldots ,t_{d-l})), \;\;\;\; k = \overline{1, d-l}$.

\end{theorem}

\begin{proof}
Each  systems (\ref{eq23}) - (\ref{eq25}) is solved locally in $V$ and any point from the obtained trajectories may serve as an initial condition for the "next" system to be locally solved in $V$ as well. The existence of local solutions is ensured by Peano theorem, which also gives an estimate of the existence intervals. 

We denote by $M = \max\left\{ |v_j|_{C(\overline{V})}, j = \overline{1, d-l} \right\}$.

Take $V_j, j = \overline{1, d-l-1}$ such that $x_0 \in V_j \subset \subset V_{j+1} \subset V$, open subsets. Let $b_1 = {\rm dist}(x^0, \partial V_1) > 0$, then we may choose $I_1 = \left[ -\displaystyle\frac{b_1}{M}, \displaystyle\frac{b_1}{M} \right]$ and the local solution of (\ref{eq23}) is obtained in $V_1$. Fix $b_2 = \min { \{ \rm dist} (x,y);$ $ x \in\partial V_1, y \in \partial V_2 \}$. Then the solution of (\ref{eq24}) exists in $I_2 = \left[ -\displaystyle\frac{b_2}{M}, \displaystyle\frac{b_2}{M} \right]$ for any initial data from $V_1$ and with the trajectory contained in $V_2$.

This argument can be iterated up to the system (\ref{eq25}) and the number of iteration steps is finite. This proves the first point.

Clearly, $y_1$ satisfies the statement on differentiability, in the origin. Then, $y_2$ satisfies it as well since $y_2 (t_1,0) = y_1 (t_1)$.  And so on, this extends step by step up to $y_{d-l}$ which is continuously differentiable in all its arguments in the origin. The formula, in the origin, follow from (\ref{eq23}) - (\ref{eq25}). For $k = d-l$, we use (\ref{eq25}) and $y_{d-l}(0,\dots,0) = y_{d-l-1}(0,\dots,0) = \dots =y_1(0) = x^0$. For $k = d-l-1$, we use the initial condition and we get

 $\displaystyle\frac{\partial y_{d-l}} {\partial t_{d-l-1}}(0, \ldots ,0) = \displaystyle\frac{\partial y_{d-l-1}} {\partial t_{d-l-1}}(0, \ldots ,0) = v_{d-l-1}(y_{d-l-1}(0,\dots,0)) = v_{d-l-1}(x^0)$. 

\noindent
This proceeds iteratively up to $k = 1$.

Due to the uniqueness property, we denote by $x^1 = y_{d-l}(t_1, \ldots ,t_{d-l})$, for some $(t_1, \ldots ,t_{d-l})$ in the existence set and by $\tilde{y}_{d-l}$, the solution of (\ref{eq23}) - (\ref{eq25}) with initial condition $x^1$. Since we are in the autonomous case with respect to all the independent variables, we have $\tilde{y}_{d-l}(s_1, \ldots ,s_{d-l}) = y_{d-l}((s_1, \ldots ,s_{d-l})+(t_1, \ldots ,t_{d-l}))$.
The proof is finished by applying the above relation in the new origin (in $x^1$).
\end{proof}

\begin{remark}\label{r22} 
\end{remark}
\vspace{-3mm}

In \cite{c11}, for $d = 3$, two iterated Hamiltonian systems are used. 
A related analysis via specific ODE's arguments together with relevant numerical examples are indicated. 
The system (\ref{eq23}) - (\ref{eq25}) is a generalization of this situation and  we underline that, as in \cite{c11}, one can approximate easily its solution, for instance with MatLab. The system (\ref{eq23}) is an usual ordinary differential system and we get its approximate solution in the discretization points of $I_1$; then the system (\ref{eq24}) is solved for each initial condition defined for the values of the parameter $t_1$ given by these discretization points in $I_1$ and so on. Finally, one obtains the approximate values of $y_{d-l}(t_1, t_2, \ldots, t_{d-l})$ on a discretization grid of $I_1 \times I_2 \times \ldots \times I_{d-l}$. In fact, all the solutions may be computed on their maximal existence interval. This can be achieved very simply and very quickly,  by standard numerical routines for ODE's.

\vspace{0.2cm}

\begin{theorem}\label{th23}
For every $k = \overline{1, l}$, $ j = \overline{1, d-l}$, we have
\begin{equation}\label{eq26}
F_k(y_j(t_1, t_2, \ldots, t_j)) = 0, \quad \forall \; (t_1, t_2, \ldots, t_j) \in I_1 \times I_2 \times \ldots \times I_j. 
\end{equation}
\end{theorem}

\begin{proof}
We notice first that, for any $k = \overline{1, l}$, we have:
$$
\displaystyle\frac{\partial}{\partial t_1}F_k(y_1(t_1)) = \nabla F_k(y_1(t_1)) \cdot v_1(y_1(t_1)) = 0, \forall t_1 \in I_1,
$$

\noindent since $v_1$ is orthogonal to $\nabla F_k$, $k = \overline{1,l}$, by (\ref{eq22}).

Moreover, $F_k(y_1(0)) = F_k(x^0) = 0$, $k = \overline{1, l}$, by (1.1). This gives (2.6) for $j = 1$. The argument follows by induction after $j$:

We assume that for $j = \overline{1, r}, r \in N$, $r \leq d-l-1$, we have (\ref{eq26}) for any $k = \overline{1, l}$ and for any $(t_1, t_2, \ldots, t_r) \in I_1 \times I_2 \times \ldots \times I_r$.

We show that this is also valid for $j = r+1$. First we remark that

\begin{equation}\label{eq27}
F_k(y_{r+1}(t_1, t_2, \ldots, t_r, 0)) = F_k(y_r(t_1, t_2, \ldots, t_r)) = 0, \; \forall \; k = \overline{1, l}, 
\end{equation}

\noindent due to the induction hypothesis.

We also notice that
\begin{equation}
\displaystyle\frac{\partial}{\partial t_{r+1}}F_k(y_{r+1}(t_1, t_2, \ldots, t_{r+1})) =\nonumber
\end{equation}

\begin{equation}\label{eq28}
= \nabla F_k(y_{r+1}(t_1, t_2, \ldots, t_{r+1})) \cdot v_{r+1}(y_{r+1}(t_1, t_2, \ldots, t_{r+1})) = 0,
\end{equation}

\begin{equation}
\forall \; (t_1, t_2, \ldots, t_{r+1}) \in I_1 \times I_2 \times \ldots \times I_{r+1}, \forall \; k = \overline{1, l}, \nonumber
\end{equation}

\noindent due to the differential equation satisfied by $y_{r+1}$ on $I_{r+1}$ and to the orthogonality relation (\ref{eq22}) satisfied by the construction of $v_{r+1}(\cdot)$. By (\ref{eq27}), (\ref{eq28}), we get (\ref{eq26}) for $j = r+1$ and this ends the proof.
\end{proof}

\vspace{0.4cm}

\noindent
Under hypothesis (\ref{eq21}), the local solution of (\ref{eq11}) is a $d-l$ dimensional manifold around $x^0$. We expect that $y_{d-l}(t_1, t_2, \ldots, t_{d-l})$ is a local parametrization of this manifold on $I_1 \times I_2 \times \ldots \times I_{d-l}$.

\begin{theorem}\label{th26}
If $F_k \in C^1(\Omega)$, $k = \overline{1, l}$, with the independence property (\ref{eq21}), and the $I_j$ are sufficiently small, $j = \overline{1, d-l}$, then the mapping
$$
y_{d-l} : I_1 \times I_2 \times \ldots \times I_{d-l} \to R^d
$$

\noindent is regular and one-to-one on its image.
\end{theorem}
\vspace{0.1cm}

\begin{proof}
\normalfont
We  get that $y_j \in C^1(I_1 \times I_2 \times \ldots \times I_j)$ by Theorem \ref{th21}.

\noindent
The matrix $B$ of partial derivatives of $y_{d-l} = (y_{d-l}^1, y_{d-l}^2, \ldots, y_{d-l}^d)$, where the superscripts denote the components of the vector $y_{d-l}$, is:

\begin{equation}\label{eq29}
B=\left(
\begin{array}{ccccc}
\displaystyle\frac{\partial y_{d-l}^1}{\partial t_1} & \displaystyle\frac{\partial y_{d-l}^1}{\partial t_2} & \ldots & \displaystyle\frac{\partial y_{d-l}^1}{\partial t_{d-l}}\\[3mm]
\ldots & \ldots & \ldots & \ldots \\[3mm]
\displaystyle\frac{\partial y_{d-l}^d}{\partial t_1} & \displaystyle\frac{\partial y_{d-l}^d}{\partial t_2} & \ldots & \displaystyle\frac{\partial y_{d-l}^d}{\partial t_{d-l}}
\end{array}
\right)
\end{equation}

\vspace{2mm}

We denote by $M_{d-l}$ the $(d-l) \times (d-l)$ matrix of the last $d-l$ rows in $B$ and we compute its determinant. Notice that the last column in $M_{d-l}$ is given by the last $d-l$ components of the vector $v_{d-l}$, that is $(0, 0, \ldots, 0, \Delta (x))^T$ due to the way we have constructed $v_{d-l}$ in (\ref{eq22}), $x$ being here the appropriate point in $V$ obtained as the value of the solution $y_{d-l}( \overline{t_1}, \overline{t_2}, \ldots, \overline{t_{d-l}})$, for some $( \overline{t_1}, \overline{t_2}, \ldots, \overline{t_{d-l}}) \in I_1 \times I_2 \times \ldots \times I_{d-l}$. We write shortly  $\Delta (y_{d-l})$ for $\Delta (x)$ with $x$ determined as above. We cut the last row and the last column in $M_{d-l}$, we denote the obtained matrix by $M_{d-l-1}$ and we have:

\begin{equation}\label{eq210}
{\rm det} M_{d-l} =  \Delta (y_{d-l}) {\rm det} M_{d-l-1}. 
\end{equation}

Taking into account the equation of $y_{d-l}$ (see (\ref{eq25})) and the fact that the components of $v_{d-l}$, from order $l+1$ to order $d-1$ are 0 (as mentioned above), the initial condition in (\ref{eq25}) gives by integration:
$$
y_{d-l}^{l+1} = y_{d-l-1}^{l+1}; \ldots; y_{d-l}^{d-1} = y_{d-l-1}^{d-1} 
$$

\noindent and they are independent of $t_{d-l}$. Therefore, we can write

\begin{equation}\label{eq211}
M_{d-l-1} = \left(
\begin{array}{ccccc} 
\displaystyle\frac{\partial y_{d-l-1}^{l+1}}{\partial t_1} & \displaystyle\frac{\partial y_{d-l-1}^{l+1}}{\partial t_2} & \ldots & \displaystyle\frac{\partial y_{d-l-1}^{l+1}}{\partial t_{d-l-1}}\\

\ldots & \ldots & \ldots & \ldots \\

\displaystyle\frac{\partial y_{d-l-1}^{d-1}}{\partial t_1} & \displaystyle\frac{\partial y_{d-l-1}^{d-1}}{\partial t_2} & \ldots & \displaystyle\frac{\partial y_{d-l-1}^{d-1}}{\partial t_{d-l-1}}
\end{array}
\right).
\end{equation}

\vspace{2mm}

Relation (\ref{eq211}) shows that in fact $M_{d-l-1}$ has a similar structure as $M_{d-l}$, associated to $y_{d-l-1}$. Using the differential system satisfied by $y_{d-l-1}$ and the structure of $v_{d-l-1}$, we see again that the last column in $M_{d-l-1}$ is of the form $(0, 0, \ldots, 0, \Delta (y_{d-l-1}))^T$ (and of length $d-l-1$). Here,  the determinant $\Delta (y_{d-l-1})$ is $\Delta (x)$ computed in the point $x = y_{d-l-1}(\overline{t_1}, \overline{t_2}, \ldots, \overline{t_{d-l-1}})$.

One can iterate the above arguments to obtain

\begin{equation}\label{eq212}
{\rm det} M_{d-l} = \Delta (y_{d-l}){\rm det} M_{d-l-1} = \Delta (y_{d-l}) \Delta (y_{d-l-1}){\rm det} M_{d-l-2} =
\end{equation}

$$
 \ldots =  \Delta (y_{d-l}) \Delta (y_{d-l-1}) \ldots \Delta (y_1) \neq 0 , 
$$

\vspace{3mm}
\noindent where the notations $M_{d-l-2}$, etc., are obvious. Relations (\ref{eq29}) - (\ref{eq212}) end the proof.
\end{proof}

\vspace{0.4cm}
\noindent
We consider now another solution choice in (\ref{eq22}). We shall use $d-l$ solutions of (\ref{eq22}) obtained by fixing the last $d-l$ components of the vector $v(x) \in R^d$ to be the rows of the identity matrix in $R^{d-l}$. The next result shows that we construct exactly the solution of the classical implicit functions theorem, which follows as a special case of our approach. 
\vspace{0.1cm}

\begin{theorem}\label{th27}
If $F_k \in C^1(\Omega)$, the last $d-l$ components of $y_{d-l}$ have the form  $(t_1 + x^0_{l+1}, t_2 + x^0_{l+2}, \ldots, t_{d-l} + x^0_{d})$, that is the first $l$ components of $y_{d-l}$ give the unique solution of the implicit system (\ref{eq11}) on $(x^0_{l+1},  x^0_{l+2}, \ldots,  x^0_{d}) + (I_1 \times I_2 \times \ldots \times I_{d-l})$ .
\end{theorem}

\begin{proof}
By inspection and induction, one can see that the last $d-l$ components of $y_j(t_1, t_2, \ldots, t_j)$ are $(t_1 + x^0_{l+1}, t_2 + x^0_{l+2}, \ldots, t_j + x^0_{l+j}, t_{j+1}+x^0_{l+j+1},
\ldots, t_{d-l}+x^0_{d})$. 

This is due to the special choice of the last components of the vectors $v_k$ in (\ref{eq22}), as rows of the identity matrix, allowing explicit integration . Then, we have just to remark that by redenoting the last $d-l$ components of  $y_{d-l}$ as $(s_{l+1}, s_{l+2}, \ldots, s_d)$, then the first $l$ components of $y_{d-l}$ are functions of $(s_{l+1}, s_{l+2}, \ldots, s_d)$, defined on $(x^0_{l+1},  x^0_{l+2}, \ldots,  x^0_{d}) + (I_1 \times I_2 \times \ldots \times I_{d-l})$, solving (\ref{eq11}) due to Theorem \ref{th23}.

The uniqueness comes from the implicit function theorem. 
\end{proof}

\begin{remark}\label{r29}
\end{remark}
\vspace{-3mm}

 We underline that, although Theorem \ref{th27} provides the classical solution of the implicit functions theorem,  a parametrization may be more advantageous in applications since it offers a more complete description of the corresponding manifold by removing the condition to obtain just functions. One can use maximal solutions of (\ref{eq23}) - (\ref{eq25}) and, in many examples, the (local) maximal solution from Theorem \ref{th26}  may give even a global description of the manifold, \cite{c11}, \cite{c17}. In applications, the choice of other solutions of (\ref{eq22}) is also possible and of interest \cite{r}, in order to improve the description of the manifold.

\begin{remark}\label{r210}
\end{remark}
\vspace{-3mm}

Beside the existence statement, Theorem \ref{th27} gives a construction recipe for the implicit functions solution and an evaluation of its existence neighborhood (via Theorem \ref{th21}), in the  system (\ref{eq11}). 

\noindent
For instance, if in the proof of Theorem \ref{th21} we take $V = B(x_0, R)$ and $V_j = B(x_0, jR(d-l)^{-1})$, then $I_j = [-R/(d-l)M, R/(d-l)M]$, for $j = 1,2,..,d-l$. This may be compared with \cite{nou1}, \cite{nou2} where other types of arguments are used.

\vspace{2mm}
\noindent
We consider now general perturbations of (\ref{eq11}) having the form

\begin{equation}\label{eq41}
F^{\lambda}_k(x_1, \ldots, x_d) = 0, \quad k = \overline{1, l}, \; \lambda \in (-1, 1) , 
\end{equation}

\noindent where $F^{\lambda}_k \in C^2(\Omega \times (-1, 1))$, $F^0_k = F_k$ and $F^{\lambda}_k(x^0) = 0$, $k = \overline{1, l}$. Hypothesis (\ref{eq21}) remains clearly valid for the perturbation as well, for $\lambda$ small.

We denote by $(S_{ \lambda} )$ the differential system similar to (\ref{eq23}) - (\ref{eq25}), associated to the perturbed implicit system (\ref{eq41}) and by $ v^{\lambda}_j$ the corresponding solutions of (\ref{eq22}), appearing in the right-hand side of $(S_{ \lambda} )$. Then, $ v^{\lambda}_j$ are in $C^1 (V_1 \times (- \lambda _0 , \lambda _0))$, under our hypotheses, for some $V_1 \in \mathcal{V}(x^0), V_1 \subset \subset V$ independent of $\lambda \in (-\lambda_0 , \lambda_0)$, for $\lambda_0$ small. The same ideas as in  Thm. \ref{th21} or Rem. \ref{r210} and the obvious property
\vspace{0.2cm}

\noindent
$M^{\lambda} = \max\left\{ |v^{\lambda}_j|_{C(\overline{V})}, j = \overline{1, d-l} \right\} \rightarrow M = \max\left\{ |v_j|_{C(\overline{V})}, j = \overline{1, d-l} \right\}$

\vspace{0.2cm}

\noindent
give the existence of the closed intervals with the origin in their interior $I_j, j = \overline{1, d-l}$, independent of $\lambda$, such that the solution of $S_{\lambda}$ is defined on $I_1 \times I_2 \times \ldots \times I_{d-l} $.

\vspace{2mm}
We denote by $y^{\lambda}_{1}(t_1), \ldots , y^{\lambda}_{d-l}(t_1, t_2, \ldots, t_{d-l})$, the unique  solution of $(S_{\lambda})$, defined in $ I_1 \times \ldots \times I_{d-l} $. By making  translations with respect to the initial conditions in each subsystem of $(S_{\lambda})$ , the initial conditions become $0$ and the differentiability properties of the solution, with respect to $ \lambda $, are a consequence of standard results on the differentiability with respect to the parameters in ODE's (since $ v^{\lambda}_j \in C^1 (V_1 \times (- \lambda _0 , \lambda _0))$) and of an inductive argument as before.
Denoting by $z^{\lambda}_{1}(t_1), \ldots , z^{\lambda}_{d-l}(t_1, t_2, \ldots, t_{d-l})$ the derivative of the above solution with respect to $ \lambda \in (-\lambda_0 , \lambda_0)$, their system in variations associated to (\ref{eq41}) and (\ref{eq11}) can be obtained by differentiation in $(S_{ \lambda} )$ with respect to $\lambda$ of the perturbations $v^\lambda_j, j = \overline{1, d-l}  $, etc.
\noindent For the case of the implicit function theorem (i.e. Theorem \ref{th27}), we obtain explicit information in algebraic form:

\begin{proposition}\label{p43}
 We have:

a) the last $d-l$ components of $z^{\lambda}_{1}(t_1), \ldots , z^{\lambda}_{d-l}(t_1, t_2, \ldots, t_{d-l})$ are null.
 
b) for any $j = 1, \ldots, d-l$ and $ \; (t_1, t_2, \ldots, t_{d-l}) \in I_1 \times \ldots \times I_{d-l} $, $z^{\lambda}_{j}(t_1, t_2, \ldots, t_{j})$ is the unique solution of:

\begin{equation}\label{eq47}
\nabla _y F_k^{\lambda} (y_j^{\lambda}) z_j^{\lambda} + \partial _{\lambda} F_k^{\lambda} (y_j^{\lambda}) = 0, \; k = 1, \ldots, l.   
\end{equation}
\end{proposition}

\begin{proof}
The first statement is a clear consequence of Theorem \ref{th27} and of the above discussion.
Since we have already established above the differentiability properties of $y_j^{\lambda}$ with respect to $\lambda$  on some given open set, one can differentiate with respect to  $\lambda$  in (\ref{eq41}) with $x_j$ replaced by $y_j^{\lambda}$, to obtain (\ref{eq47}). Notice that the solution of the linear system (\ref{eq47}) is unique due to (\ref{eq21}) and to a).
\end{proof}

\begin{remark}\label{r44}
\end{remark}
\vspace{-3mm}

One can obtain for $z^{\lambda}_{1}(t_1), \ldots , z^{\lambda}_{d-l}(t_1, t_2, \ldots, t_{d-l})$ the relation (\ref{eq47}) even for implicit parametrizations as in Theorem \ref{th26}, but point $a)$ is not valid and (\ref{eq47}) is not uniquely determining (without supplementary information)
 $z^{\lambda}_{1}(t_1),\ldots,
 z^{\lambda}_{d-l}(t_1, t_2, \ldots, t_{d-l})$ . 
To obtain the necessary supplementary information, one has to use directly the differential systems (\ref{eq23})-(\ref{eq25}) and to compute the corresponding system in variations. 

\vspace{2.5mm}
\noindent
Consider now, as an example, the special case of perturbations of the form

\vspace{-4mm}
\begin{equation}\label{eq42}
F_j(x_1, \ldots, x_d) + \lambda h_j(x_1, \ldots, x_d) = 0, \quad j = \overline{1, l}, \; \lambda \in (-1, 1) , 
\end{equation}
\noindent
where $ h_j \in C^2 (\Omega), \ h_j(x^0) = 0$. 

If, moreover, $l = 1$ and the equation $ F(x_1, \ldots, x_d) = 0$, $F \in C^2(\Omega )$, together with the associated initial condition, represents the boundary of a subdomain in $\Omega$ (where $F<0$, for instance) then the geometric perturbation defined by (\ref{eq42}) may be very complex, including topological and boundary perturbations \cite{c8}, \cite{hp}, \cite{sz}. Computing the equation in variations as in Proposition \ref{p43}, the perturbations (\ref{eq42}) generate a directional derivative in the implicit system (\ref{eq11}). Consequently, by the above geometric interpretation, we may define, for $l = 1$, a new type of geometric directional derivative of domains. This is more general than the speed method or the topological derivatives \cite{c8} and has applications in shape optimization, fixed domain methods, see \cite{pt}, \cite{ti}.

\vspace{0.1cm}

\section{Generalized solutions}

In this section, we discuss the problem (\ref{eq11}) for $F_j \in C^1(\Omega), j = \overline{1,l}$, in the absence of the hypothesis (\ref{eq21}) i.e. all determinants of maximal order $l$ may be null in $x^0$. We remark that there is $\{x^n\} \subset \Omega$, such that:

\begin{equation}\label{eq31}
x^n \to x^0, \quad {\rm rank} J(x^n) = l, \; n \in N,
\end{equation}

\noindent where $J(x^n)$ denotes the Jacobian matrix of $F_1, F_2, \ldots, F_l \in C^1(\Omega)$, in $x^n$.

Notice that in case (\ref{eq31}) is not fulfilled, it means that rank $J(x) < l$ in $x \in W$, where $W$ is a neighborhood of $x^0$. Then $F_1, F_2, \ldots, F_l$ are not functionally independent in $W$ and the problem (\ref{eq11}) can be reformulated by using less functionals \cite{c12}, \cite{c14}. That is (\ref{eq31}) is in fact always valid, except for not well formulated problems, including redundant equations. One may classify the systems of type (\ref{eq11}), from this point of view, in well-posed and ill-posed systems. Notice as well that (\ref{eq31}) is fulfilled if (\ref{eq21}) holds, i.e. (\ref{eq31}) is the generalization of (\ref{eq21}), valid for all well-posed implicit systems.

Due to (\ref{eq31}), in each $x^n$, one can use the results of the previous section for the system

\begin{equation}\label{eq32}
F_j(x) - F_j(x^n) = 0, \; j = \overline{1, l}, \; x \in \overline{\Omega}, 
\end{equation}

\noindent where we can find locally the solution of (\ref{eq32}) around $x^n$, in a neighborhood depending on $n$.

From (\ref{eq31}), we also have $F_j(x^n) \to F_j(x^0) = 0$, for $ n \to \infty$, $j = \overline{1, l}$, since $F_j \in C^1(\Omega)$.

We denote by $T_n$ the closure in $\overline{\Omega}$ of the manifold defined by (\ref{eq32}). It  is compact and connected. We also have that $\{T_n\}$ are uniformly bounded since $\Omega$ is bounded and, on a subsequence denoted by $\alpha$, we get

\begin{equation}\label{eq33}
T_n \to T_{\alpha}, \; n \to \infty, 
\end{equation}

\noindent in the Hausdorff-Pompeiu metric \cite{c8}, \cite{c7}, where $T_{\alpha}$ is some compact connected subset in $R^d$.

\vspace{0.2cm}

\begin{definition}\label{def31}
$T =  \mathop{\bigcup}\limits_{\alpha} T_{\alpha}$ is the (local) generalized solution of (\ref{eq11}) in $x^0$. The union is taken for all the sequences and subsequences satisfying (\ref{eq31}), (\ref{eq33}).
\end{definition}
\vspace{0.1cm}

This notion was introduced in \cite{c17} and further discussed in \cite{c11}, in dimension two and three, by exploiting continuity properties with respect to data in Hamiltonian systems. The present treatment in arbitrary dimension is based on general convergence properties and allows a relaxation of the regularity conditions.

The above definition covers all critical or non critical cases.
\noindent
See Remark \ref{r34} as well. For instance, if in (\ref{eq11}) we have just one equation and $x^0$ is an isolated extremum for the respective function, then the generalized solution is just $\{x^0\}$. If the respective function is identically zero in the open set $O \subset \Omega$ and $x^0$ is on the boundary of $O$, then (\ref{eq31}) is satisfied and the generalized solution is the boundary of $O$ or some subset of it - see Proposition \ref{p33} and Example \ref{ex35} below.  A complete description of the level sets (even of positive Lebesgue measure)  may be obtained in this way via the generalized solutions. The generalized solution is not a manifold and may be not a compact subset (for instance, if $\Omega$ is unbounded), but it is connected, \cite{c8}, Appendix 3. The approximating generalized solution, i.e. $\bigcup  T_{n_0}$ (for  some "big" $n_0$ in (\ref{eq33}) and for several choices of the approximating sequences of $x_0$ in (\ref{eq31})), may be not connected. One can easily approximate the generalized solutions, by the techniques from Section 2 applied to the corresponding terms from the sequence $ \{x_n\} $ close enough to $x^0$. Due to the properties of the Hausdorff-Pompeiu distance, the approximation is uniform in the space variables. If not enough sequences are taken into account, it is possible to obtain (locally) just a subset of $T$. For instance, in the equation $x^2 - y^2 = 0$, around the origin, with one approximating sequence  $(x_n,y_n) \to (0,0)$, such that $|y_n|<x_n$, just some part of the solution is generated at the limit. Taking into account a supplementary sequence such that $x_n<-|y_n|$  the whole solution is obtained (locally) by Definition \ref{def31}. An algorithm for the approximation of the generalized solution is discussed in \cite {r}, including many relevant examples.

\vspace{0.2cm}
\noindent
Let $M \subset \overline{\Omega}$ denote the connected component of the solution of (\ref{eq11}), containing the critical point $x^0$. If $intM $ is nonvoid, then it does not contain $x^0$, due to (\ref{eq31}), that is $x^0 \in \partial M$.

\vspace{1mm}
\begin{proposition}\label{p33}
We have:  $x^0 \in T_{\alpha} \subset T \subset \partial M_{x_0}, \forall \alpha$,
where $ \partial M_{x_0}$ is the connected component of $\partial M$ containing $x_0$. In particular

\begin{equation}\label{eq34}
F_j(x) = 0, \; j = \overline{1, l}, \; \forall \; x \in T. 
\end{equation}
\vspace{0.1cm}
\end{proposition}

\begin{proof}
By (\ref{eq31}), we have $x^n \in T_n, \; \forall \; n$ and we get $x^0 \in T_{\alpha}$ by the definition of the Hausdorff-Pompeiu convergence. The next inclusion follows by Definition \ref{def31}. 

The same argument gives that, for any $x \in T$, then $x \in T_{\beta}$ for some subsequence $\beta$, and there are $\lambda_n \in T_n$ (here $T_n$ is the subsequence convergent to $T_{\beta}$)  such that $\lambda_n \to x$ for $n \to \infty$. By (\ref{eq32}), we see that $F_j(\lambda_n) = F_j(x^n) \to F_j(x^0) = 0$, $j = \overline{1, l}$, on a subsequence. Then, by continuity, $F_j(\lambda_n) \to F_j(x) = 0$ as claimed and (\ref{eq34}) is proved. 

Consequently, $T_{\alpha} \subset M, \forall \alpha$. If $int M$ is nonvoid, then it is formed just of points not satisfying (\ref{eq21}) since $\nabla F_j$ are null.  Then $\{x^n\}$ are disjoint from  $\overline{int M}$ ($M$ is not necessarily a Caratheodory set and may be distinct from $\overline{int M}$) and, consequently, $T_{\alpha} \subset \partial M$. By Prop.A3.2 in \cite{c8}, each $T_{\alpha}$ is connected and  contains $x_0$, by the above argument. If $\partial M$ has more connected components, then it yields $T_{\alpha} \subset \partial M_{x_0}, \forall {\alpha}$. Definition \ref{def31} ends the proof.
\end{proof}

\vspace{0.1cm}

\begin{remark}\label{r34} 
\end{remark}
\vspace{-3mm}
\noindent
If $x^0$ is a regular point, i.e. (\ref{eq21}) is satisfied, then we denote by $S$ the manifold giving the (local) solution of (\ref{eq11})  around $x^0$. Then $S$ coincides with the generalized solution around $x^0$.

\noindent
In Definition \ref{def31}, we may choose $x^n \to x^0$, $x^n \in S$ and the uniqueness property from the implicit functions theorem gives (for this choice) that $T_n = S$ locally, for $n$ big enough. This choice of $\{ x^n \} $ satisfies (\ref{eq31}) since $J(x^n) \to J(x^0)$, so $ x^n  $ satisfies (\ref{eq21}) for $n$ big enough. We see that in the classical case, one obtains $T = S$ (locally), that is Definition \ref{def31} gives indeed a generalization of the classical local solution of the implicit functions theorem.

\vspace{0.1cm}

\begin{example}\label{ex35}

In $R^2$, take $d = 2, l = 1$ and 

\begin{equation}\label{eq35}
f(x_1, x_2) = \left\{
\begin{array}{lll}
x_1^2(x_2^2 - x_1^2)^2 & if \; x_1 < 0, &|x_2| \leq |x_1|\\
0 & otherwise.
\end{array}\right.
\end{equation}

\end{example}

Clearly $f$ is in $C^1(R^2)$ and $\nabla f(x_1, x_2) = 0$, on the second line of (\ref{eq35}). Take $x^0 = (0, 0)$ and $x^n \to x^0$, $x^n = (x_1^n, x_2^n)$, $x_1^n < 0$, $|x_2^n| < |x_1^n|$.

In such points $x^n$, one can use Theorem \ref{th26} and (\ref{eq22}), together with the relations (\ref{eq23}) - (\ref{eq25}), give the Hamiltonian system (in dimension two, iterated systems are not necessary):

\begin{equation}\label{eq36}
\begin{array}{cc}
x_1'(t) = -4x_1^2x_2(x_2^2 - x_1^2),\\

x_2'(t) = 2x_1(x_2^2 - x_1^2)(x_2^2 - 3x_1^2),\\

(x_1(0), x_2(0)) = x^n.
\end{array} 
\end{equation}

Here, we have chosen $(-f_{x_2}, f_{x_1})$ as the solution of (\ref{eq22}).

\begin{figure}[htbp]
	\centering
		\includegraphics[width=0.7\textwidth]{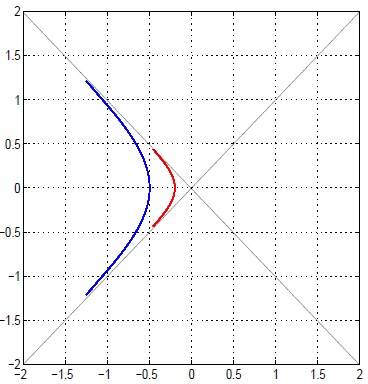}
		\caption{}
	\label{fig1}
\end{figure}

\vspace{3mm}
In Figure \ref{fig1}, we represent the solution $T_n$ of (\ref{eq36}) obtained with MatLab, for $x^n = (-\displaystyle\frac{1}{n}, 0), n = 2,5$. The generalized solution of the implicit function problem (\ref{eq11}) corresponding to (\ref{eq35}) is given by $T = \{ (x_1, x_2) \in R^2; x_1 = \pm \; x_2, x_1 \leq 0 \}$, the boundary of the critical set of $f(\cdot, \cdot)$, to which $x^0$ belongs.

The generalized solution contains the essential information about the solution set of (\ref{eq11}), since it gives its boundary (and in Proposition \ref{p33} the inclusion becomes equality, in this example). 

 If we define

$$
f_1(x_1,x_2) = x_1^2 [(x_1^2 + x_2^2 -1)_+ ]^2 
$$

\noindent
and $x^0 = (0,1)$, then $\partial M$ is connected and the corresponding generalized solution is $\partial M$ without the lower half of the unit circle. The inclusion in Proposition \ref{p33} is strict and $M$ is not Caratheodory,  in this case. This is also related to the local character of the construction from Section 2. See Ex. 2 in \cite{c11} as well.

We continue now with a partial converse of Proposition \ref{p33} that shows that the notion of generalized solution is a strict extension of the classical notion of solution. 

\begin{proposition}\label{p36}
Let $x^0$ be the unique critical point of (\ref{eq11}) in the interior of the closed ball $B(x^0)$. Then, $T = M$ in $B(x^0)$.
\end{proposition}

\vspace{-3mm}
\begin{proof}
Due to Proposition \ref{p33}, we have just to prove $M \subset T$.

Let $A$ be a connected component of $M - \{x^0\}$. It is open in the relative topology of $M \cap \overline{B(x^0)}$ since all the points except $x^0$ are regular and the implicit functions theorem can be applied. It is also maximal in the sense that it cannot be strictly extended in $\overline{B(x^0)}$. Notice that in the relative topology of $M$, we have $\partial A \subset \partial B(x^0) \cup \{x^0 \}$, by the implicit function theorem. Consequently, $\overline A \subset A \cup \{x^0 \}$ since $A$ is maximal and the part of $\partial A$ contained in $\partial B(x^0)$ is also contained in $A$.

We have $x^0  \in \overline A$. Otherwise, by the above relation, it yields $A = \overline A$, that is $A$ is both closed and open in $M$ and this contradicts $M$ connected.

One can consider a sequence ${x^n} \in A, x^n \to x^0$ and the associated manifolds $T_n$. Notice that  $ A = T_n $ by the implicit functions theorem. It follows that $A \cup \{ x^0 \} = lim \; \overline{T}_n ,\; \; A \cup \{ x^0 \} \subset T$. As $A$ is an arbitrary component of $M - \{x_0\}$, we get the conclusion and finish the proof. 
\end{proof}

\vspace{2mm}

\begin{proposition}\label{p37}
Let $F_j \in C^1(\Omega)$, $j = \overline{1, l}$ and $x^n \to x^0$, $x^n, x^0 \in \Omega$. Denote by $\widetilde{T}_n, \widetilde{T}_0$ the generalized solutions of (\ref{eq11}) contained in the bounded domain $\Omega$, corresponding to the initial conditions $x^n$, respectively $x^0$. Then

\vspace{-0.5cm}

\begin{equation}\label{eq37}
\mathop{\lim\sup}\limits_{n \to \infty} \widetilde{T}_n \subset \widetilde{T}_0. 
\end{equation}
\end{proposition}

\begin{proof}
Let $\widehat{x}_{n_k} \in \widetilde{T}_{n_k}$, $\widehat{x}_{n_k} \to \widehat{x}$, where $n_k \to \infty$ is some subsequence. We show that $\widehat{x} \in \widetilde{T}_0$.

By Definition \ref{def31}, there is $\widetilde{x}_{n_k} \in \Omega$, such that (\ref{eq21}) is satisfied in $\widetilde{x}_{n_k}$ and $|\widetilde{x}_{n_k} - x^{n_k}| < \displaystyle\frac{1}{n_k}$ (here, we also use the characterization of the Hausdorff-Pompeiu limit) and there are $y_{n_k} \in T_{\tilde{x}_{n_k}}$ such that $|y_{n_k} - \widehat{x}_{n_k}| < \displaystyle\frac{1}{n_k}$. Consequently, $y_{n_k} \to \widehat{x}$ fo $n_k \to \infty$. Here, $T_{\tilde{x}_{n_k}}$ is the solution of (\ref{eq32}) corresponding to $\tilde{x}_{n_k}$. By using the sequences $\widetilde{x}_{n_k} \to x^0$ and $y_{n_k} \in T_{\tilde{x}_{n_k}}$, $y_{n_k} \to \widehat{x}$, we see that $\widehat{x} \in T_0$ due to Definition \ref{eq31} and the proof is finished.
\end{proof}
\vspace{-.1cm}

\begin{example}\label{ex38}
Let $x_n \to 0, x_n < 0$ be some strictly increasing sequence and $g_n : R^2 \to R$ be given by

\vspace{-0.3cm}

\begin{equation}
g_n(x, y) = \left\{ \begin{array}{l}
c_n[(x - x_n)^2 + y^2 - \displaystyle\frac{1}{4}\min\{|x_{n+1} - x_n|^2; |x_n - x_{n-1}|^2 \}]^2 , if\nonumber\\
 \; |x - x_n|^2 + |y|^2 \leq \displaystyle\frac{1}{4}\min\{|x_{n+1} - x_n|^2; |x_n - x_{n-1}|^2\}\nonumber\\
0 \quad otherwise,\nonumber
\end{array}\right.
\end{equation}
\end{example}

\noindent where $c_n > 0$ is some "big" constant. We consider the function $F : R^2 \to R$ by

\begin{equation}\label{eq38}
F(x, y) = f(x, y) + \mathop{\sum}\limits_{n = 1}^{\infty} g_n(x, y), 
\end{equation}

\noindent where $f$ is given in (\ref{eq35}). Clearly $F$ is in $C^1(R^2)$ and $(x_n, 0)$ are local maximum points of $F$ if $c_n$ are "big". The sum in (\ref{eq38}) has always just maximum two non zero terms due to the form of the supp $g_n$.

Take the sequence $x^n = (x_n, 0) \to (0, 0) = x^0$. Then, the implicit equations $F(x, y) = F(x_n, 0)$ have the unique solution $x^n = (x_n, 0)$ in a neighbourhood of $x^n$ and $T_{x^n} = (x_n, 0)$. In the point $(0,0)$, we have $T_0$ as in Example \ref{ex35}. We see in this example that the inclusion in (\ref{eq37}) may be strict.



\section{Reduced gradients in nonlinear programming}

In constrained optimization, projected gradient methods are a classical tool, but their application may be hindered by the difficulty to effectively compute projections on the admissible set, Ciarlet \cite{ci}. Based on the results from the previous sections, we use here the reduction approach to eliminate, totally or partially, the constraints (and the Lagrange multipliers), that allows optimality conditions in a more effective way, decreasing the dimension. Local and global algorithms and numerical examples are also discussed, under weak assumptions. The elimination of certain unknowns has advantages at computational level.

In the recent papers \cite{ssb}, \cite{mcb}, dimensional reduction is obtained via new relaxation procedures associated to implicit functions. Our approach is certainly different and ensures good numerical results. In the case of polynomial and semi-algebraic optimization, \cite{la} Thm.6.5, Thm.7.5, in the setting of global optimization, a stronger constraint qualification is used.

\vspace{2mm}

We consider now the classical  minimization problem with equality constraints:

\vspace{2mm}
$(P) \;\;\;\;\;\;\;\; Min\{ g(x_1, \dots ,x_d)\}$

\vspace{1mm}
\noindent
subject to (\ref{eq11}).  It is known that by Theorem \ref{th27}  we can replace it (around $x^0$) by the unconstrained problem for $(t_1, t_2, \dots , t_{d-l}) \in (I_1 \times I_2 \times \ldots \times I_{d-l})$: 

\vspace{2mm}
$(P_1) \;\; Min\{ g(y_{d-l}^1,y_{d-l}^2, \dots , y_{d-l}^l,  t_1 + x^0_{l+1}, t_2 + x^0_{l+2}, \ldots, t_{d-l} + x^0_{d})\}$,

\vspace{2mm}
\noindent
where $(y_{d-l}^1,y_{d-l}^2, \dots , y_{d-l}^l,  t_1 + x^0_{l+1}, t_2 + x^0_{l+2}, \ldots, t_{d-l} + x^0_{d})$ are the components of $y_{d-l}$, the solution of
(\ref{eq23})-(\ref{eq25}), corresponding to this case. 
This methodology can be extended to the case of implicit parametrizations.  

By Theorem \ref{th27}, Theorem \ref{th21} and the chain rule, one easily obtains the (known) first order optimality conditions in the Fermat form, involving the tangential gradient to the constraints manifold:

\begin{proposition}\label{c48}
If $x^0$ is a local solution of $(P)$ satisfying that g and $F_i , i= \overline{1,l}$,  are in $C^1(R^d)$ and (\ref{eq21}) holds, then we have:

\begin{equation}\label{eq48}
\nabla g(x^0).v_j(x^0) = 0 \;\; j = \overline{1, d-l}.
\end{equation}

\end{proposition}

\noindent
In fact, this is equivalent with the classical Lagrange multipliers rule, since under (\ref{eq48}),  $\nabla g(x^0)$ is in the normal space, which has the basis given by $\nabla F_i(x^0), i = \overline{1, l}$.

\vspace{2mm}
In this non convex setting, we introduce the following algorithm of projected gradient type, based on the use of the tangential gradient:

\begin{algorithm}\label{b}

1) choose $n=0$, $\delta > 0$ (a tolerance parameter) and denote by $t^n = (t_1^n,\dots,t^n_{d-l})$ such that $y_{d-l}(t_1^n,\dots,t^n_{d-l}) = x^n$ in (\ref{eq23})-(\ref{eq25}).

2) compute $\rho^{n+1} \in [0, \alpha_n]$ via the line search:

\vspace{1mm}
$Min \hspace{1mm} g[y_{d-l}(t^n - \rho [\nabla g(x^n) . v_j(x^n)]_{j = \overline{1, d-l}})]$.
\vspace{1mm}

3) set:

$x^{n+1} = y_{d-l}(t^n - \rho^{n+1} [\nabla g(x^n) . v_j(x^n)]_{j = \overline{1, d-l}})$,

$t^{n+1} = t^n - \rho^{n+1} [\nabla g(x^n) . v_j(x^n)]_{j = \overline{1, d-l}}$.

4) If $|g(x^{n}) - g(x^{n+1})| < \delta$, then STOP! Otherwise n:=n+1 and GO TO Step2).

\end{algorithm}

\begin{remark}\label{r155}
\end{remark}
\vspace{-3mm}
The algorithm works practically in $V$, where the system (\ref{eq23})-(\ref{eq25}) is defined and the parameter $\alpha_n$ in the line search with limited minimization rule has to be chosen "small", such that we remain in $V$ and the system (\ref{eq23})-(\ref{eq25}) can be solved around $t^n$, in $Step \; 2)$. In $Step \; 3)$ we perform the "projection" on the constraints manifold $M \subset \Omega$. The points $x^n$ generated by this algorithm are always admissible for $(P)$. No convexity properties are assumed. The definition of $(P_1)$ uses the implicit function Theorem \ref{th27} which is appropriate for optimality conditions, while for the Algorithm \ref{b} the general implicit parametrization method has to be taken into account. The same is valid for the subsequent problem ${(Q_1)}$ and the related results.

\vspace{2mm}

In this algorithm, $\Omega$ is a bounded domain, $g \in C^1(\Omega)$ is bounded from below and the constraints are as in (\ref{eq11}) with hypothesis (\ref{eq21}) satisfied in $x^0$. We denote by $G(t) = g(y_{d-l}(t))$, defined in a neighborhood of the origin in $R^{d-l}$ and of class $C^1$ due to (\ref{eq23})-(\ref{eq25}) and Theorem \ref{th21}.  The sequence $\{g(x^n) = G(t^n)\}$ is non increasing and convergent in this general setting, ensuring the convergence of the algorithm. The sequence $\{x^n\}$ is bounded. Moreover, we have $\nabla G(t^n) = [\nabla g(x^n) . v_j(x^n)]_{j = \overline{1, d-l}}$ by Theorem \ref{th21} and the Algorithm \ref{b} is in fact a transcription of the classical gradient method for the unconstrained problem $(P_1)$. One can discuss other (very rich) variants of such local algorithms with their convergence (to stationary points, in general), under supplementary hypotheses if necessary, Bertsekas \cite{ber}, Patriksson \cite{pa}. The new point in Algorithm \ref{b} is that one can effectively compute the "projection" $y_{d-l}$.

\vspace{2mm}

We discuss now the general case of both equality and inequality constraints:

\vspace{2mm}
$(Q) \;\;\;\;\;\;\;\; Min\{ g(x_1, \dots ,x_d)\}$

\vspace{1mm}
\noindent
subject to (\ref{eq11}) and to

\vspace{-6mm}
\begin{equation}\label{eq50}
G_j (x) \leq 0  \;\; j = \overline{1,m},
\end{equation}

\noindent
where $g, F_i, G_j$ are in $C^1(R^d)$. The Mangasarian-Fromovitz condition in this case consists of (\ref{eq21}) and there is $d \in R^d$ such that

\begin{equation}\label{eq51}
\nabla F_i (x^0)d = 0, \; i = \overline{1,l}, \;\;  \nabla G_j (x^0)d < 0 , \;\;  j \in I(x^0),
\end{equation}

\noindent
with $I(x^0)$ being the set of indices of active inequality  constraints in $x^0$. See \cite{bs}, \$ 2.3.4 or \cite{Cl}, \$ 6 for excellent presentations. The necessary and sufficient metric regularity condition from \cite{tz} cannot be used here due to the lack of convexity.

The reduced problem is again obtained via Theorem \ref{th27}:

\vspace{2mm}
$(Q_1) \;\; Min\{ g(y_{d-l}^1,y_{d-l}^2, \dots , y_{d-l}^l,  t_1 + x^0_{l+1}, t_2 + x^0_{l+2}, \ldots, t_{d-l} + x^0_{d})\}$,

\vspace{2mm}
\noindent
subject to the constraints (\ref{eq50}), in the "reduced" form:

\vspace{-4mm}
\begin{equation}\label{eq52}
G_j (y_{d-l}^1,y_{d-l}^2, \dots , y_{d-l}^l,  t_1 + x^0_{l+1}, \ldots, t_{d-l} + x^0_{d}) \leq 0  \;\; j = \overline{1,m},
\end{equation}

\begin{lemma}\label{l50}
The minimization problem $(Q_1)$ satisfies the Mangasarian-Fromovitz condition in the origin of $R^{d-l}$.

\end{lemma}

\begin{proof}
\noindent
By the first part in (\ref{eq51}), we see that $d$ is in the tangent space to the manifold (\ref{eq11}) since $\nabla F_i (x^0), i = \overline{1,l}$ is a basis in the normal space to the manifold given ($\ref{eq11}$), under hypothesis ($\ref{eq21}$).  Then $d = \mathop{\sum}\limits_{s = 1}^{d-l}\alpha_s v_s$ with $\alpha_s$ some scalars, since $v_s, s = \overline{1,d-l}$, gives a base in the tangent space.

\noindent
By the second part in (\ref{eq51}) we get $\mathop{\sum}\limits_{s = 1}^{d-l}\alpha_s \nabla G_j (x^0) v_s < 0$. Using the derivation formula from Theorem \ref{th21}, this may be rewritten as $\mathop{\sum}\limits_{s = 1}^{d-l}\alpha_s \displaystyle\frac{\partial}{\partial t_s} g_j (0,0, \dots, 0) < 0$, where

$g_j (t_1, \dots, t_{d-l}) = G_j(y_{d-l}^1,y_{d-l}^2, \dots , y_{d-l}^l,  t_1 + x^0_{l+1}, \ldots, t_{d-l} + x^0_{d})$.

\vspace{2mm}
\noindent
is the composed mapping. This shows that the Mangasarian-Fromovitz hypothesis is satisfied in the origin of $R^{d-l}$
with the vector $(\alpha_1, \dots, \alpha_{d-l})$.

\end{proof}

\vspace{1.5mm}
If $x^0$ is a local solution of $(Q)$, by Lemma \ref{l50}, one can apply the classical KKT theorem, \cite{ci}, to the problem $(Q_1)$ in the origin of $R^{d-l}$ that becomes a local solution for $(Q_1)$. Using again the derivation formula, we get:
 
\begin{theorem}\label{th50}
Let $x^0$ be a local minimum for $(Q)$. Then, there are $ \beta_j \geq 0, j = \overline{1, m}$ such that

$0 = \nabla g(x^0).v_s(x^0) + \mathop{\sum}\limits_{j = 1}^{m} \beta_j \nabla G_j (x^0).v_s(x^0), s = \overline{1,d-l}$,

$0 = \beta_j G_j (x^0), j = \overline{1,m}$.
\end{theorem}

\begin{remark}\label{r55}
\end{remark}
\vspace{-2mm}
\noindent
This is a simplified version of the KKT conditions since it eliminates the Lagrange multipliers for the equality constraints.
It is possible to eliminate completely the Lagrange multipliers:
if $x^0$ is a local solution of problem $(Q)$, then one can remove  the inactive inequality constraints at $x^0$.
This is a consequence of the remark that the inequality constraints that are not active at $x^0$ define a neighborhood of $x^0$.  The minimum property of $x^0$  is preserved in this neighborhood, just under the equality constraints supplemented by the active constraints rewritten as equalities. Under the independence condition for all these constraints, one can write optimality conditions as in the Proposition \ref{c48}.

\vspace{2mm}
We relax now the hypotheses in the problem $(Q)$ and we describe a direct minimization algorithm of global type. It looks for the solution in a maximal neighborhood of $x^0$, corresponding to the maximal solutions of the subsystems in (\ref{eq23}) - (\ref{eq25}) (the maximal existence intervals may depend on the respective initial conditions). See Remark \ref{r29} and \cite{c11}, \cite{c17}.

We assume in the sequel that $g$ and $G_j, j = \overline{1,m}$, are just in $C(R^d)$ and $F_i,  i = \overline{1,l}$, are in $C^1(R^d)$ and satisfy condition (\ref{eq21}) in $x^0$. This last condition can be removed in fact, working with generalized solutions, according to the subsequent Remark \ref{r60}. Notice that $x^0$ is here just an admissible point for  $(Q)$ and not a local minimum as in Theorem \ref{th50}. We can also add the abstract constraint $x \in D$, some given subset in $R^d$, such that $x^0 \in D$.

The main observation is that in solving numerically (\ref{eq23}) - (\ref{eq25}), now using the variant corresponding to Theorem \ref{th26}, we obtain automatically a discretization of the manifold defined by (\ref{eq11}), in a maximal neighborhood of $x^0$, as explained above. Let us denote by $n$ the discretization parameter. For instance, $1/n$ can characterize the size of the discretization for the parameters $t_1, \dots, t_{d-l}$, $n$ or may be linked to the length of the intervals where the maximal solution is computed, etc.
We denote by $C_n$ the set of all these discretized points that, moreover, satisfy all the constraints (the inequality and the other restrictions have to be just checked). They give the approximating admissible set and we formulate the  algorithm:

\begin{algorithm}\label{a}

1) choose $n = 1$, the discretization step $1/n$ and the 

solution intervals $I_1^n, \dots, I_{d-l}^n$, the tolerance parameter $\delta$.

2) compute the discrete set of admissible points $C_n$, starting from $x^0$, 

via (\ref{eq23}) - (\ref{eq25}) and by testing the validity of  (\ref{eq50}) and $D$.

3) find in $C_n$ the approximating minimum of $(Q)$, denoted by   $x^n$. 

4) test if the solution is satisfactory by $|g(x_n) - g(x_{n-1} | \leq \delta$.

5) If YES, then STOP. If NO, then $n := n+1$ and GO TO step 1).

\end{algorithm}

In step 4) other tests (on the solutions, on the gradients, etc.) may be used. The approximating minimum $x^n \in C_n$ may be not unique and the Algorithm \ref{a} finds all all of them. One can adapt the convergence test to such situations. 

\begin{theorem}\label{th55}
The algorithm is convergent as $n \rightarrow \infty$.

\end{theorem}

This is a consequence of the density of $\bigcup C_n$ in the admissible set, according to Theorem \ref{th26}.

\begin{remark}\label{r65}
\end{remark}
\vspace{-3mm}

The set defined by the equality constraints may have several connected components. See Example \ref{e65}.  Starting from $x^0$, Algorithm \ref{a} will minimize just on the component that contains $x^0$. Initial guesses from all the admissible components are necessary if we want to minimize on all of them.

\begin{remark}\label{r60}
\end{remark}
\vspace{-3mm}

If condition (\ref{eq21}) is not fulfilled , then one can use the generalized solution of (\ref{eq11}) as explained in Section 3 (see Proposition \ref{p36}) , since the Hausdorff-Pompeiu distance ensures the uniform convergence of approximating points. The computed minimum may satisfy (\ref{eq11}) or the minimum property with some small error tolerance and the convergence property with respect to the discretization parameters is ensured. An algorithm for the computation of the generalized solution, with relevant examples is studied in \cite{r}.

\vspace{2mm}
\noindent
Finally, we indicate some illustrative  numerical examples and compare our results with other approximation methods, from MatLab or \cite{ssb}.

\begin{example}\label{e60}
\end{example}

\vspace{-3mm}
We consider first a minimization problem on the torus in $R^3$, with radii 2 respectively 1,  defined implicitly by $F = 0$, and with initial point $(x_0,y_0,z_0) =(\sqrt{5},2,0)$:

\vspace{-0.2cm}

\begin{eqnarray*}
&&min\{xyz\}\\
&&F(x,y,z)=(x^2+y^2+z^2+3)^2-16(x^2+y^2)
\end{eqnarray*}

The obtained results are given below, compared with the application of the fmincon routine of MatLab:

\vspace{-0.2cm}

\begin{eqnarray*}
&&min=-2,7154\\
&&x_{min}=1,7841; y_{min}=1,8199; z_{min}=-0,8363\\
&&fmincon: min=-2,7153; x_{min}=1,802; y_{min}=1,802; z_{min}=-0.836
\end{eqnarray*}

Using other starting points like $(1,0,0)$ or $(3,0,0)$ is not allowed by MatLab that finds no other admissible solutions in these cases, while our approach works.

\begin{example}\label{e65}
\end{example}

\vspace{-2mm}
Now, we consider two equality restrictions, given by $F$ and $ P$, that represent a torus intersected with a paraboloid, see Fig.2 and Fig.3. Two initial points are taken into account since the intersection has two components.

\begin{eqnarray*}
&&min\{x^3+5y-7sin z\}\\
&&P(x,y,z)=\displaystyle\frac{2\sqrt{3}}{3}x-y^2-z^2\\
&&(x_0,y_0,z_0)=(\sqrt{3},1,1); \;(x_0,y_0,z_0)=(\sqrt{3},-1,1)
\end{eqnarray*}

The numerical results and a comparison with MatLab routine fmincon is indicated below:

\vspace{-0.2cm}

\begin{eqnarray*}
&&(\sqrt 3,-1, 1):      minimal\; value = 0.498975897823261\\
&&solution:\; (1.06688905550184, -0.814925789648031, 0.753631933331335)\\
&&(\sqrt 3,1, 1):     minimal\; value = -7.65929313197537 \\
&&solution:\; (1.10697710321061, 0.817093948780941, 0.781479124977557)
\end{eqnarray*}

In the second case fmincon stops after 42 iterations with the message that constraints are not satisfied within the tolerance. In the first case, fmincon finds basically the same solution.

\vspace{4mm}
\begin{remark}\label{r80}
\end{remark}

\vspace{-3mm}
In \cite{ssb}, an example in $R^6$, with three equality constraints, is discussed. Reworking it via Algorithm \ref{a}, starting from the two points indicated there on p.451,  we obtain the new points

$(0.5631, -3.2581, 0.51593, 0.4692, 1.4635, 3.589)$,

$(0.56166,-3.3154, 0.50897, 0.5047, 1.4365, 3.6777)$ 

\noindent
with the cost values $343,7695$, respectively $383,7265$. This improves the quoted experiment and can be directly checked. It does not contradict \cite{ssb} since our algorithm needs no bounds on the independent variables and extends the search domain, which is an advantage from the point of view of global optimization. The necessary working time, on a medium performance laptop, is several minutes. More details on the experiment and some high dimensional numerical examples are indicated in \cite{R2}.

\begin{figure}
	\centering
		\includegraphics[width=0.7\textwidth]{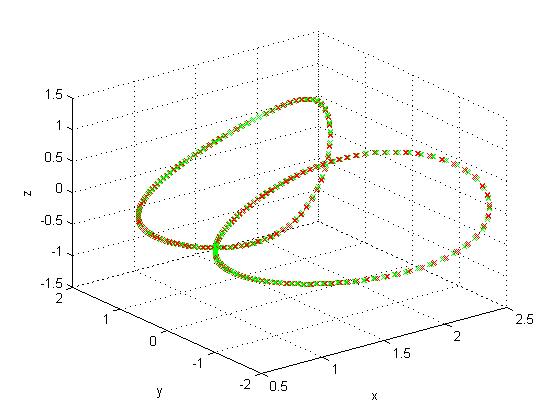}
		\caption{The admissible set}
	\label{fig2}
\end{figure}

\begin{figure}
	\centering
		\includegraphics[width=0.7\textwidth]{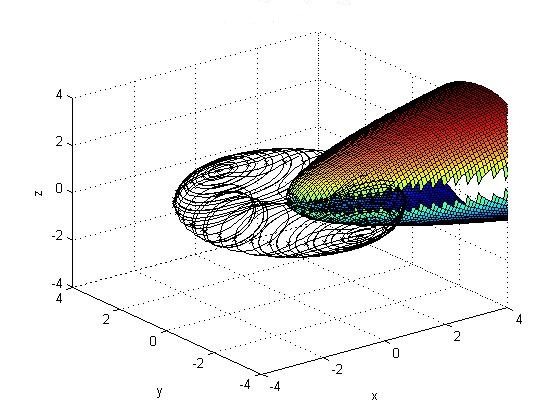}
		\caption{The geometry}
	\label{fig3}
\end{figure}


\newpage
\begin{thebibliography}{}

\bibitem{ber} D. P. Bertsekas, Nonlinear Programming (3rd Edition), Athena Scientific, Nashua, MASS (2016).

\bibitem{bs} J.F. Bonnans and A. Shapiro, Perturbation analysis of optimization problems, Springer Verlag, New York (2000).

\bibitem{bd} F. Bouchut, L. Desvillets, On two-dimensional Hamiltonian transport equations with continuous coefficients, Diff. Int. Eqns., vol.14, no.8, pp. 1015-1024 (2001).

\bibitem{a} G.J. Butler and H.I. Freedman, Further critical cases of the scalar implicit function theorem, Aequationes Math. 8, pp. 203–211 (1972).

\bibitem{nou1} H.C. Chang, W. He and N. Prabhu, The analytic domain in the implicit function theorem, JIPAM, Vol. 4, Iss. 1, Article 12 (2003).

\bibitem{ci} Ph. Ciarlet, Introduction to numerical linear algebra and optimization, Cambridge Univ. Press, New York (1989).

\bibitem{Cl} F.H. Clarke, Optimization and nonsmooth analysis, John Wiley \& Sons, New York (1983).

\bibitem{b} E. Coddington and N. Levinson, Theory of ordinary differential equations, McGraw-Hill, New York (1955).

\bibitem{dz} M. Delfour and J.-P. Zolesio, Shapes and Geometry, SIAM, Philadelphia (2001).

\bibitem{dpl} R.J. DiPerna, P.L. Lions, Ordinary differential equations, transport theory and Sobolev spaces, Invent. Math. 98, pp.511-547 (1989).

\bibitem{c2} K. Dobiasova, Parametrizing implicit curves, WDS'08 Proceedings of Contributed Papers, MATHFYZPRESS, pp. 19-22, Prague (2008).

\bibitem{c3} A.L. Dontchev and R.T. Rockafellar, Implicit functions and solution mappings, Springer, New York (2009).


\bibitem{c4} Xiao-Shan Gao, Conversion between implicit and parametric reprezentations of algebraic varieties, Mathematical mechanization and applications, Academic Press, pp. 253-271, San Diego (2000).


\bibitem{hp} A. Henrot, M. Pierre, Variation et optimization de formes: une analyse geometrique, Springer Verlag, Berlin (2005).

\bibitem{c6} S.G. Krantz and H.R. Parks,The implicit function theorem, Birkh\" auser, Boston (2002).

\bibitem{kt} C. Kublik and R. Tsai,  Integration over curves and surfaces defined by the closest point mapping, preprint ArXiv 1504.05478v4 (2015).

\bibitem{c7} C. Kuratowski, Introduction to set theory and topology, Pergamon Press, Oxford (1962).

\bibitem{la} J.B. Lasserre, An introduction to polynomial and semi-algebraic optimization, Cambridge University Press, Cambridge (2015).

\bibitem{c} S. Lefschetz, Differential equations: geometric theory, Interscience, New York (1957).



\bibitem{mcb} A. Mitsos, B. Chachuat, P.I. Barton, McCormick-based relaxations of algorithms, SIAM J. Optim. 20(2), pp.573-601 (2009).

\bibitem{c8} 
\newblock P. Neittaanm\" aki, J. Sprekels, D. Tiba, Optimization of elliptic systems. Theory and applications, Springer, New York (2006).

\bibitem{pt} P. Neittaanmaki, D. Tiba, Fixed domain approaches in shape optimization problems, Inverse Problems, vol.28, p.1-35, (2012)
doi:10.1088/0266-5611/28/9/093001

\bibitem{c11} M.R. Nicolai and D. Tiba, Implicit functions and parametrizations in dimension three: generalized solutions, DCDS - A vol. 35, no.6, pp.2701 - 2710 (2015). doi:10.3934/dcds.2015.35.2701

\bibitem{r} 
\newblock M.R. Nicolai, An algorithm for solving implicit systems in the critical case, Ann. Acad. Rom. Sci. Ser. Math. Appl. Vol. 7, no. 2, pp.310 - 322, (2015).

\bibitem{R2} 
\newblock M.R. Nicolai, High dimensional applications of implicit parametrizations in nonlinear programming, Ann. Acad. Rom. Sci. Ser. Math. Appl., Vol. 8, no.1, pp.44 -55, (2016).

\bibitem{c12} M. Nicolescu, N. Dinculeanu and S. Marcus, Analiz\u a Matematic\u a, vol. I, Ed. 4, Ed. Didactic\u a \c si Pedagogic\u a, Bucure\c sti (1971).

\bibitem{pa} M. Patriksson, Nonlinear Programming and Variational Inequality Problems: A Unified Approach, Springer (2013).


\bibitem{nou2} Phan Phien, Some quantitative results on Lipschitz inverse and implicit functions theorems, arXiv: 1204.4916v2 (2012).

\bibitem{c14} W. Rudin, Principles of mathematical analysis, Second Edition, McGraw-Hill, New York (1964).

\bibitem{c15} J. Schicho, Rational parametrizations of algebraic surfaces, Thesis, J. Kepler Univ. Linz (1995).

\bibitem{sz} J. Sokolowski, J.-P. Zolesio, Introduction to shape optimization. Shape sensitivity analysis, Springer Verlag, Berlin (1992).

\bibitem{ssb} M.D. Stuber, J.K. Scott, P.I. Barton, Convex and concave relaxations of implicit functions, Optimization methods and software, 30(3), pp.424-460 (2015).

\bibitem{c16} J.A. Thorpe, Elementary topics in differential geometry, Springer Verlag, New York (1979).

\bibitem{c17} D. Tiba, The implicit functions theorem and implicit parametrizations, Ann. Acad. Rom. Sci. Ser. Math. Appl. 5, no. 1-2, pp. 193-208, (2013). http://www.mathematics-and-its-applications.com

\bibitem{ti} D.Tiba, Boundary Observation in Shape Optimization, in "`New trends in differential equations, control theory, and optimization"', V. Barbu, C. Lefter, I. Vrabie (Eds.), World Scientific Publishing, Singapore, pp.301 - 314 (2016).

\bibitem{tz} D. Tiba, C. Zalinescu, On the necessity of some constraint qualification conditions in convex programming, J.Convex Anal. vol.11, no.1, pp.95-110 (2004).

\bibitem{c18} D. Wang, Irreducible decomposition of algebraic varieties via characteristic set method and Gr\" obner basis method, CAGD 9, pp. 471-484, (1992).

\bibitem{c19} H. Yang, B. J\" uttler, L. Gonzales-Vega, An evolution-based approach for approximate parametrization of implicitly defined curves by polynomial parametric spline curves,
\newblock Math. Comp. Sci. 4, no. 4, pp. 463-479 (2010)

\end{thebibliography}


\end{document}